\documentclass[graybox]{svmult}

\usepackage{mathptmx}       
\usepackage{helvet}         
\usepackage{courier}        
\usepackage{makeidx}         
\usepackage{graphicx}        
\usepackage{multicol}        
\usepackage[bottom]{footmisc}
\usepackage{amssymb}  
\usepackage{amsmath}

\makeindex             

\newcommand{\R}{\mathbb{R}}
\newcommand{\Z}{\mathbb{Z}}

\newcommand{\take}{\!\setminus\!}

\begin{document}

\title*{Some Comments on the Inverse Problem \\ of Pure Point Diffraction}
\author{Venta Terauds and Michael Baake}

\institute{Fakult\"{a}t f\"{u}r Mathematik, Universit\"{a}t Bielefeld, 
Postfach 100131, 33501 Bielefeld, Germany\newline \and 
 Venta Terauds: terauds@math.uni-bielefeld.de
\newline \and Michael Baake:   mbaake@math.uni-bielefeld.de}

%
%
\maketitle

\abstract*{In a recent paper \cite{TeBa-lenzmoo}, Lenz and Moody
  presented a method for constructing families of real solutions to
  the inverse problem for a given pure point diffraction
  measure. Applying their technique and discussing some possible
  extensions, we present, in a non-technical manner, some examples of
  homometric structures.}

\abstract{In a recent paper \cite{TeBa-lenzmoo}, Lenz and Moody
  presented a method for constructing families of real solutions to
  the inverse problem for a given pure point diffraction
  measure. Applying their technique and discussing some possible
  extensions, we present, in a non-technical manner, some examples of
  homometric structures.}

\section{Introduction}

Kinematic diffraction is concerned with the Fourier transform
$\widehat{\gamma}$ of the autocorrelation $\gamma$ of a given
structure, the latter described by a measure $\omega$, which is
usually assumed to be translation bounded; see \cite{TeBa-baagrm12}
for a recent summary of the state of affairs. Of particular relevance
(for crystals and quasicrystals, say) are systems with pure point (or
pure Bragg) diffraction, i.e., those where $\widehat{\gamma}$ is a
pure point measure.

In considering the inverse problem, namely the problem of determining
which structure or structures could have produced a given diffraction,
one is naturally led to the concept of \emph{homometry}, where
structures that give rise to the same diffraction measure are said to
be homometric.  Via the autocorrelation measure, homometry is
well-defined for certain classes of measures, and, accordingly, also
for objects such as point sets and tilings (an overview of these
concepts is given in \cite[Ch.~9.6]{TeBa-book}). Various methods have
been used (see \cite{TeBa-baagrm06, TeBa-book, TeBa-baagrm08,
  TeBa-grumoo} and references therein) to construct different objects
(mainly Dirac combs) that are homometric. However, a way of finding
\emph{all} objects with a given diffraction, which we shall refer to
as the \emph{diffraction solution class} of a given diffraction
measure, has long remained elusive.

In \cite{TeBa-lenzmoo}, Lenz and Moody present a method for abstractly
parametrising the real solution class of a given pure point
diffraction measure. Their approach is based on the Halmos-von Neumann
theorem in conjunction with concepts from the theory of (stochastic)
point processes.  The objects constructed via their method are in many
cases measures, but as we shall see, even for a very simple periodic
diffraction, one may construct objects with that diffraction that
generically fail to be measures.  Thus many open questions remain, in
particular whether (or in what sense) non-measure solutions have a
reasonable physical interpretation.

The mathematical formalism behind the construction method in
\cite{TeBa-lenzmoo} is based on a theory of (generalised) point processes
and is quite formidable. In essence, it justifies the use of the
lower path of the Wiener diagram 
\begin{center}
\begin{picture}(90,60)(0,8)
 \put(10,19){$\widehat{\omega}$}
 \put(10,60){$\omega$}
 \put(70,61){$\gamma$}
 \put(70,19){$\widehat{\gamma}$}
 \put(21,62){${\vector(1,0){46}}$}
 \put(38,63){$\stackrel{\circledast}{ }$}
 \put(21,22){${\vector(1,0){46}}$}
 \put(38,23){$\stackrel{|\cdot|^2}{ }$}
 \put(14,56){${\vector(0,-1){25}}$}
 \put(4,40){$\stackrel{\mathcal{F}}{ }$}
 \put(73,56){${\vector(0,-1){25}}$}
 \put(75,40){$\stackrel{\mathcal{F}}{ }$}
\end{picture}
\end{center}\vspace*{-4mm}
in the reverse direction, where $\mathcal{F}$ denotes Fourier
transform and $\circledast$ the volume-averaged (or Eberlein)
convolution (so that $\gamma = \omega \circledast \widetilde{\omega}$;
compare \cite{TeBa-hof97,TeBa-baagrm12}). In particular, for a
diffraction measure of the form $\widehat{\gamma} = \sum_{k\in
  L}|A(k)|^2 \delta_k$, with $L$ a countable set, one may (formally,
but consistently) invert the left as well as the bottom arrow of the
diagram (the latter interpreted as $|\cdot|^2$ being applied to the
coefficients or `intensities' $I(k) = \lvert A(k) \rvert^{2}$
individually) to construct a measure (or at least a tempered
distribution) $\omega$ with diffraction $\widehat{\gamma}$.

It is well known that for certain measures (for example the Dirac comb
of a lattice or, more generally, of a crystallographic structure
\cite{TeBa-hof95,TeBa-hof97,TeBa-franz}), one may proceed via the
lower route in the forward, and thus also the reverse, direction in
this way. However, the results of Lenz and Moody apply in a much more
general situation.

The purpose of this brief contribution is to present, without the
point process formalism of \cite{TeBa-lenzmoo}, some examples of object
classes that display the same diffraction.  Unless stated otherwise,
all measures presented below are measures on the real line $\R$, which
means that we illustrate everything with one-dimensional examples.

\section{A simple diffraction measure with simple origins}

A \emph{diffraction measure}, that is, a measure that represents the
diffraction of some (physical) structure, must be real, positive and
inversion symmetric.  By \emph{pure point}, we mean that the measure
can be written as $\widehat{\gamma} = \sum_{k\in L} I(k) \delta_k$,
with a countable set $L$ (which may be finite) and locally summable
intensities $I(k)>0$.  Let us begin then, with the `simplest
possible' pure point diffraction measure.

\begin{example}
  Let $\widehat{\gamma} = \delta_0$.  Proceeding backwards through the
  Wiener diagram via the `bottom' route, we gain $\widehat{\omega} =
  A(0)\delta_0$, where, applying the method of \cite{TeBa-lenzmoo}, we
  must have $A(0)=1$.  Since $\widehat{\lambda} = \delta_0$ and
  $\widehat{\delta_0} = \lambda$, where $\lambda$ is Lebesgue measure,
  we have $\omega = \lambda$ and the real diffraction solution class
  of $\delta_0$ is $\{\lambda, -\lambda\}$ (as it is clear that
  $-\omega$ is homometric to $\omega$). Of course, observing that
  $\lambda\circledast\widetilde{\lambda} = \lambda$ gives us the same
  thing via the top route in the Wiener diagram.

  It is not hard to deduce a bit more here:\ For $u\in \mathbb{S}^1$,
  the unit circle, we have $u\lambda\circledast \widetilde{u\lambda} =
  u\overline{u}\lambda = \lambda$, and so, via the top route, we see
  that $\{u\lambda\,|\, u\in\mathbb{S}^1\}$ is contained in the
  (complex) diffraction solution class of $\delta_0$. Are there
  further measures with diffraction $\delta_0$? Well, any measure of
  the form $\omega + \mu$, with $\mu$ a finite measure on $\R$, has
  the same diffraction as $\omega$, as adding a finite measure to
  $\omega$ does not change the autocorrelation. This is well-known
  \cite{TeBa-hof95,TeBa-book} and a `trivial' degree of freedom; in
  the framework of \cite{TeBa-lenzmoo}, the point process for $\omega
  + \mu$ is the same as that for $\omega$.
\end{example}

In this example, we have a good idea of what the complex diffraction
solution class is, and are certain that the real diffraction solution
class contains only measures.  In the next section, we shall see that
the real diffraction solution class of a nice, periodic measure like
$\delta_{\Z}$ contains both measures and non-measures. In fact, as
will be shown in \cite{TeBa-ter}, the only pure point diffraction measures
whose real solution class consists solely of measures are those
supported on a finite set of points.

\section{A lattice diffraction measure with all kinds of  origins}

Let us now consider measures with diffraction $\widehat{\gamma} =
\delta^{}_{\Z} = \sum_{k\in\Z} \delta^{}_k$. According to the Wiener
diagram, objects $\omega$ with the diffraction $\delta^{}_{\Z}$ must
have the (possibly formal) Fourier transform $\widehat{\omega} =
\sum_{k\in\Z} A(k) \delta^{}_k$, with $|A(k)| = 1$ for all $k$. In the
setting of \cite{TeBa-lenzmoo}, we have the further conditions that $A(0) =
1$ and $A(-k) = \overline{A(k)}$ for all $k$, and one constructs
different objects with diffraction $\delta^{}_{\Z}$ simply by choosing
different sets of compliant coefficients $\{A(k) \mid k \in \Z \}$. An
interpretation of $\omega$ and $\widehat{\omega}$ might need the
theory of tempered distributions and their relations with measures.

\begin{example}
  Choosing $A(k) = 1$ for all $k$, one gains $\omega = \delta_{\Z}$,
  as follows from the Poisson summation formula $\widehat{\delta^{}_{\Z}}
  = \delta^{}_{\Z}$; compare \cite{TeBa-cordo,TeBa-book}.
\end{example}

\begin{example}
  By splitting the set $\Z$ into subsets $n\Z, n\Z+1, \ldots,
  n\Z+(n-1)$, and choosing coefficients appropriately, one may
  construct an $n$-periodic measure with diffraction
  $\delta^{}_{\Z}$. For example, to construct a $4$-periodic measure,
  let
\begin{equation*} 
   A(k)\, = \, \left\{
	\begin{array}{ll}
	    1 	\,,		& \quad k\in 4\Z \, , \\
	    \alpha \,, 		& \quad k\in 4\Z+1 \, , \\
	    e	\,, 		& \quad k\in 4\Z + 2 \, , \\
	    \overline{\alpha}\,,& \quad k\in 4\Z + 3 \, ,
          \end{array}        \right.               
\end{equation*}
where $e = \pm 1$ and $\alpha\in\mathbb{S}^1$. Note that this
exhausts all possibilities for compliant coefficients for this
four-way splitting.  We have
\begin{equation*} 
   \widehat{\omega^{\pm}_{\alpha}}\, = \, (\delta^{}_0 + \alpha\delta^{}_1 +
   \overline{\alpha}\delta^{}_{-1} \pm \delta^{}_2) \ast \delta^{}_{4\Z} \,.
\end{equation*}
Then, writing $\alpha = e^{2\pi \I t_\alpha}$ and defining 
$\chi^{}_s(t):= e^{2\pi \I st}$ for $t,s\in\R$, we obtain
\begin{align*}
   \omega^{\pm}_{\alpha} &\, = \, \tfrac{1}{4}(\chi^{}_0 + \alpha\chi^{}_1 + 
   \overline{\alpha}\chi^{}_{-1} \pm \chi^{}_2)\delta_{\frac{1}{4}\Z} \\
			&\, = \, \tfrac{1}{4}\sum_{k\in\Z} 
   \bigl(1 + 2\cos(2\pi(t_\alpha + \tfrac{k}{4})) \pm \cos(\pi k)\bigr) 
   \delta_{\frac{k}{4}} \,,
\end{align*}
where we have used the the Poisson summation formula and the identity
$\widehat{\mu \ast \nu} = \widehat{\mu}\widehat{\nu}$.  That is, we
have a measure of the form $\omega^{\pm}_{\alpha} = (a\delta^{}_0 +
b\delta_{\frac{1}{4}} + c\delta_{\frac{1}{2}} +
d\delta_{\frac{3}{4}})\ast\delta^{}_{\Z}$.  For example, choosing
$\alpha \in \{0,\I,1,-\I\}$, we get $\omega^{\pm}_{\alpha}$ according
to the following table.

\begin{center}
\begin{small}
\begin{tabular}{|r|l|l|}
\hline
\parbox[c][5mm][c]{0cm}{}$\omega_{\alpha}^{\pm}\;\;\;$ 
        &\; $e=1$ &\; $e=-1$ \\
\hline
\parbox[c][6mm][c]{0cm}{} $t_{\alpha} = 0$ \;& \;$\delta^{}_{\Z}$ 
      & \;$\frac{1}{2}(\delta^{}_0 + \delta_{\frac{1}{4}} - 
        \delta_{\frac{1}{2}} + \delta_{\frac{3}{4}})\ast\delta^{}_{\Z}$ \;\\
\hline
\parbox[c][6mm][c]{0cm}{}\;$t_{\alpha} = \frac{1}{4}$\,\, &  
         \;$\frac{1}{2}(\delta^{}_0 - 
       \delta_{\frac{1}{4}} + \delta_{\frac{1}{2}} + 
       \delta_{\frac{3}{4}})\ast\delta^{}_{\Z}$  \;
      & \;$\delta_{\frac{3}{4}}\ast\delta^{}_{\Z}$ \\
\hline
\parbox[c][6mm][c]{0cm}{}$t_{\alpha} = \frac{1}{2}$\,\, &  
     \;$\delta_{\frac{1}{2}}\ast\delta^{}_{\Z}$
      & \;$\frac{1}{2}(-\delta^{}_0 + \delta_{\frac{1}{4}} + 
      \delta_{\frac{1}{2}} + \delta_{\frac{3}{4}})\ast\delta^{}_{\Z}$ \;\\
\hline
\parbox[c][6mm][c]{0cm}{}$t_{\alpha} = \frac{3}{4}$\,\, & 
     \;$\frac{1}{2}(\delta^{}_0 + \delta_{\frac{1}{4}} + 
     \delta_{\frac{1}{2}} - \delta_{\frac{3}{4}})\ast\delta^{}_{\Z}$ 
      & \;$\delta_{\frac{1}{4}}\ast\delta^{}_{\Z}$ \\
\hline
\end{tabular}\;.
\end{small}
\end{center}
\medskip 

\noindent
It is easy to verify (along the top route of the Wiener diagram) that
these measures do indeed all have autocorrelation
$\gamma=\delta^{}_{\Z}$ and thus diffraction $\widehat{\gamma} =
\delta^{}_{\Z}$. All elements constructed are pure point measures;
compare \cite{TeBa-grumoo} for similar examples.
\end{example}

\begin{example}
  We now construct a measure that is not itself pure point, but has
  diffraction $\delta^{}_{\Z}$. Take $A(0) = 1$ and $A(k) = -1$ for all
  $k\in\Z\take\{0\}$. Then we have
\begin{equation*} 
    \widehat{\omega} \;=\; \delta^{}_0 \, - 
    \!\!\!\!\sum_{k\in\Z\,\take\{0\}} \!\!\!\delta^{}_k 
    \;=\; 2\delta^{}_0 - \delta^{}_{\Z}\,,
\end{equation*}
and see that $\omega = 2\lambda - \delta^{}_{\Z}$ is a member of the
(real) diffraction solution class of $\delta^{}_{\Z}$.
\end{example}

Recalling that $\lambda\circledast\delta^{}_{\Z} = \lambda$, one can
easily verify the above example and construct similar ones via a
little trial and error and the top route of the Wiener diagram. For
example, via a quick calculation with the Eberlein convolution, one
can see that the diffraction of $\delta^{}_{\Z} -(1+\I)\lambda$ is
also $\delta^{}_{\Z}$, while $\delta^{}_{\Z} - \lambda$ has
diffraction $\delta^{}_{\Z\,\take\{0\}}$. Such a trial and error
method, however, would probably not lead to the following
more elaborate system.

\begin{example}
  To construct something a little different, we use an aperiodic set
  suggested by the period doubling sequence, which is limit-periodic;
  compare \cite{TeBa-book,TeBa-franz}. From the set
\begin{equation*}
    \Lambda \, = \, \bigcup_{n\geq 0} (2.4^n \Z + (4^n - 1))\, 
\end{equation*}
of left endpoints of `$a$'s in this sequence \cite{TeBa-baamoo04}, we
form the symmetric set
\begin{equation*} 
  \varDelta \;:=\; \Lambda \cup (-\Lambda) \;=\; 2\Z \,
  \mathbin{\dot{\cup}} \mathbin{\dot{\bigcup_{n\geq 1}}} \varDelta_n\,, 
\end{equation*}
where $\varDelta_n := \bigl(2.4^n \Z + (4^n - 1) \bigr) \cup
\bigl(2.4^n \Z + (1-4^n)\bigr)$ for $n\geq 1$ and $\dot{\cup}$ denotes
the disjoint union of sets.  Then, we define $A(k) = 1$ for $k\in
\varDelta$ and $A(k) = -1$ for $k\notin \varDelta$.  Due to the
symmetry of $\varDelta$, these coefficients satisfy the conjugacy
condition given above. We get
\begin{equation*} 
     \widehat{\omega} \;=\, \sum_{k\in \varDelta} \delta^{}_k - 
    \!\sum_{k\in \Z\,\take\varDelta} \!\!\delta^{}_k \;=\; 2 
    \delta^{}_{\varDelta} - \delta^{}_{\Z} \,,
\end{equation*}
and thus $ \omega = 2 \widehat{\delta^{}_{\varDelta}} -
\delta^{}_{\Z}$.  Now, we have
\begin{equation*}
     \delta^{}_{L} \, = \; 
     \delta^{}_{2\Z} + \sum_{n\geq 1} \delta^{}_{\varDelta_n} 
	   \, = \; \delta^{}_{2\Z} + \sum_{n\geq 1} \delta^{}_{2\cdot 4^n \Z} 
              \ast (\delta^{}_{4^n - 1} + \delta^{}_{1-4^n})\,,
\end{equation*}
and hence may use the Poisson summation formula to derive the formal
expression
\begin{align*}
\widehat{\delta^{}_{\varDelta}} 	&\, = \, \tfrac{1}{2}\delta_{\frac{\Z}{2}} +
          \sum_{n\geq 1} \frac{\overline{\chi}_{4^n - 1} + 
         \overline{\chi}_{1-4^n}}{2.4^n} \delta_{\frac{\Z}{2.4^n}} \\
			&\, = \, \tfrac{1}{2}\delta_{\frac{\Z}{2}} + 
         \sum_{n\geq 1} \frac{\cos(2\pi(4^n - 1)k)}{4^n}
          \delta_{\frac{\Z}{2.4^n}}\,.
\end{align*}
As $\delta^{}_\varDelta$ is a translation bounded measure,
$\widehat{\delta^{}_{\varDelta}}$ is a tempered distribution. It is
not, however, a measure, as it is easy to find compact sets
$K\subseteq\R$ (for example, take $K= [0,\frac{1}{4}]$) such that
$|\widehat{\delta^{}_{\varDelta}}|(K)$ is infinite. So, we have a
tempered distribution, $\omega$, that is not a measure, but
nevertheless has diffraction $\delta^{}_{\Z}$.
\end{example}

Of course, one can construct many such `non-measures' (in
\cite{TeBa-ter}, via a theorem of Cordoba \cite{TeBa-cordo}, it is
shown that the homometry class of $\delta^{}_{\Z}$ contains
uncountably many such objects), but this does not shed much light on
the physical relevance of such constructions. A little more insight
may be gained by noting that (in this case at least) our constructed
distribution is the limit (in the weak-$\ast$ topology on the space
$\mathcal{S}^{\prime} (\R)$ of tempered distributions) of a sequence
of measures over $\R$.

For $\epsilon > 0$, define  
\begin{equation*}
    \rho^{}_{\varepsilon} \, := \, 
    \tfrac{1}{2}\delta_{\frac{\Z}{2}} + \sum_{n\geq 1} \frac{\cos(2\pi(4^n
  - 1)k)}{(4+\varepsilon)^n} \delta_{\frac{\Z}{2.4^n}}\,.
\end{equation*}
A short calculation reveals that, for $\varepsilon>0$,
$|\rho^{}_{\varepsilon}|(K)$ is finite for all compact sets $K\subseteq
\R$, so that $\rho^{}_{\varepsilon}$ is indeed a measure (it is even
translation-bounded). Thus $\omega^{}_{\varepsilon}:= 2\rho^{}_{\varepsilon}
- \delta^{}_{\Z}$ is also a measure. Moreover, for all Schwartz functions
$g\in\mathcal{S} (\R)$, one has $\omega^{}_{\varepsilon}(g) \to
\omega(g)$ as $\varepsilon \to 0^+$. This is a standard approach in
Fourier analysis to enforce convergence of the series, which is
sometimes referred to as `regularisation' in physics. Such objects can
still be given a reasonable physical meaning.

\section{Further remarks}

The method of \cite{TeBa-lenzmoo} may only be applied to diffraction
measures $\widehat{\gamma}$ that are `backward transformable',
meaning that the (inverse) Fourier transform, $\gamma$ (that is, the
autocorrelation), is again a measure. Further contemplation of the
previous example, however, makes this condition seem a little too
restrictive. Using the scheme, we constructed the object $\omega
= 2\widehat{\delta^{}_{\varDelta}} - \delta^{}_{\Z}$, which is not a
measure, but is the weak-$\ast$ limit (in $\mathcal{S}' (\R)$) of
measures. The object $\widehat{\delta^{}_{\varDelta}}$ is also a
non-measure weak-$\ast$ limit of measures, but does not have the good
fortune, as $\omega$ does, to have a measure-valued autocorrelation,
so it is not covered by the scheme.

Presuming that one may extend the method of \cite{TeBa-lenzmoo} to
admit this case (or, in other words, that one may proceed via the
lower path in the Wiener diagram to calculate the diffraction), one
has $\omega = \widehat{\delta^{}_{\varDelta}}$, with $\widehat{\omega}
= \delta^{}_{\varDelta}$, and thus diffraction $\widehat{\gamma} =
\delta^{}_{\varDelta}$. The measure $\delta^{}_{\varDelta}$ is
positive, inversion symmetric and translation bounded, so has almost
all of the properties that one expects from a diffraction measure.

The natural next step is to understand the classes of objects for
which one may define an autocorrelation (and thus a diffraction). The
framework of \cite{TeBa-lenzmoo} is applicable in the more general
setting of a locally compact Abelian group. However, if one considers
only objects in $\R^d$, these initial examples suggest that
consideration of tempered distributions that are the weak-$\ast$
limits of measures may be a good place to begin.

\begin{acknowledgement}
  This work was supported by the German Research Council (DFG), via the
  CRC 701, and by the RCM$^2$, at the University of Bielefeld.
\end{acknowledgement}

\end{document}